\pdfoutput=1
\documentclass[11pt]{article}

\usepackage[a4paper,hmargin=1in,vmargin=1.25in]{geometry}
\usepackage[utf8]{inputenc}
\usepackage[toc,page]{appendix}
\usepackage{multicol}
\usepackage{amsmath,amsfonts,amssymb,amsthm,url}
\usepackage{tikz}
\usepackage{longtable}
\usepackage{booktabs,authblk}
\usepackage{multirow}
\usepackage{threeparttable}
\usetikzlibrary{shapes,arrows}
\usetikzlibrary{backgrounds}

\newcommand*{\lon}{
       \mskip1mu
        \relax
        {;}
        \mskip1mu
        \relax
}

\usepackage{caption}
\usepackage{subcaption}

\usepackage{longtable}
\usepackage{mathtools}

\def\Real{\mathbf{R}} 

\def\kos{K}

\def\sss{\smallsetminus}

\def\cip{\texttt{CI}}

\theoremstyle{plain}
\newtheorem{proposition}{Proposition}[section]
\newtheorem{theorem}[proposition]{Theorem}

\theoremstyle{definition}

\title{Tensor Product of Polymatroids and Common Information}
\author{Carles Padr\'o}
\affil{Universitat Polit\`ecnica de Catalunya, Barcelona, Spain}

\begin{document}

\maketitle

\begin{abstract}
A new connection between two
different necessary conditions for a polymatroid to be 
linearly representable is presented.
Specifically, we prove that  
the existence of a tensor product with the uniform matroid
of rank two on three elements implies 
the existence of a  common information 
extension for every pair of subsets of the ground set.
\end{abstract}

\section{Introduction}

This work is motivated by a question
posed in~\cite[Section~5]{BGILMS25}
about the connections
between two necessary conditions for a 
polymatroid to be linearly representable.
Namely, the  existence of common information 
extensions~\cite{BFP25,DFZ09,HRSV00}
and the existence of tensor products
with regular matroids~\cite{BGILMS25,LsV81,Lov77}.

Another  necessary condition is 
given by Ingleton's inequality~\cite{Ing71}.
Its original proof~\cite{Ing71}
deals with intersections of 
two and three vector subspaces. 
A simpler proof using only one intersection of 
two subspaces was presented in~\cite{HRSV00}.
As a consequence, the inequality holds
for every polymatroid
that admits common information extensions.
Subsequently, other linear rank inequalities 
were found~\cite{DFZ09} by using the 
common information property.

An alternative proof for the Ingleton's inequality
has been recently presented~\cite{BGILMS25}.
Specifically, a polymatroid satisfies the
Ingleton's inequality if
it admits a tensor product with the
uniform matroid $U_{2,3}$,
which implies that the inequality holds 
for every linearly representable polymatroid. 

Some questions arise from that proof,
as the one posed in~\cite{BGILMS25}: 
which are the connections
between common information and tensor products?
Or are those two necessary conditions essentially different?
Is it possible to find 
new linear rank inequalities by using tensor products?

Some answers to those questions are presented in this work.
Mainly, we prove that the existence of a tensor
product with $U_{2,3}$ implies the
existence of a common information extension
for every pair of subsets of the ground set
(Theorem~\ref{st:main}).
That result was found after noticing that the 
values of rank function of the tensor product of 
a linearly representable polymatroid
with $U_{2,3}$ depend on the dimensions of
intersections of three vector 
subspaces (Proposition~\ref{st:reptens}).
The bounds on the rank function 
that are derived from that fact 
apply to every tensor product with $U_{2,3}$
(Proposition~\ref{st:gentens}).
Those bounds provide in particular
a simpler and more transparent proof for the fact that
the Ingleton's inequality is satisfied
by every polymatroid that admits 
a tensor product with $U_{2,3}$
(Proposition~\ref{st:newIng}).

\section{Preliminaries}

A compact notation for set unions is used
and we avoid the curly brackets for singletons.
that is, we write $XY$ for $X \cup Y$ 
and $X y$ for $X \cup \{y\}$.
Similarly, $X \sss x$ denotes 
the set difference $X \sss \{x\}$. 
Given $A \subseteq E$ and $i \in \{1,2,3\}$, 
we notate $A^i$  for the subset
$A \times \{i\}$
of $E \times \{1,2,3\}$.
For example, $A^1 B^2 (CD)^3$
denotes the set 
$(A \times \{1\}) \cup 
(B \times \{2\}) \cup
((C\cup D)\times \{3\})$.
For a \emph{set function}
$f \colon 2^E \to \Real$
and subsets $X,Y,Z \subseteq E$, we write
\[
f(Y \lon Z| X) =
f(XY) + f(XZ) - f(XYZ) - f(Z)
\]
and, in particular,
\(
f(Y \lon Z) = 
f(Y \lon Z | \emptyset) =
f(Y) + f(Z) - f(YZ)
\)
and
\(
f(Y | X) = f(Y \lon Y | Z) = 
f(X Y) - f(X)
\). 

A \emph{polymatroid}
is a pair $(E,f)$ formed by a 
finite \emph{ground set} $E$ and a
\emph{rank function} $f$ that is a 
monotone and submodular set function on $E$
with $f(\emptyset) = 0$.
If $f$ is integer-valued and 
$f(X) \le |X|$ for every $X \subseteq E$, 
then $(E,f)$ is a \emph{matroid}.
We are going to use the
following result 
from~\cite[Proposition~44.1]{Sch03}.

\begin{proposition}
\label{st:polymchar1}
For a set function $f$ on $E$,  
the following properties are equivalent.
\begin{enumerate}
\item
$f$ is monotone and submodular.
\item
\(
f(Y\lon Z|X) \ge 0
\)
for every $X,Y,Z \subseteq E$.
\item
$f(y\lon z |X) \ge 0$ for every 
$X \subseteq E$ and $y,z \in E \sss X$.
\end{enumerate}
\end{proposition}

If $(E,f)$ is a polymatroid and
$Z$ is finite set with 
$Z \cap E = \emptyset$, 
a polymatroid $(EZ, g)$ such that
$g(X) = f(X)$ for every $X \subseteq E$ is 
called an \emph{extension} of $(E,f)$.
The same symbol will be used most of the times
for the rank functions of a polymatroid and 
its extensions.

Given a collection $(V_x)_{x \in E}$
of vector subspaces of a 
vector space $V$ over a field $\kos$, 
the set function $f$ on $E$ defined by 
$f(X) = \dim(\sum_{x \in X} V_x)$ is the 
rank function of a \emph{polymatroid} on $E$.
In that situation, the collection $(V_x)_{x \in E}$
is a \emph{$\kos$-linear representation} 
of the polymatroid $(E,f)$.
In that situation, the polymatroid is 
said to be \emph{$\kos$-linearly representable}, or
simply \emph{$\kos$-linear}.

\section{Necessary Conditions for Linearity}

By using linear algebra, one can prove in different 
ways~\cite{DFZ09,HRSV00,Ing71} that, for every 
four subsets $A,B,C,D$ of the ground set
of a linear polymatroid,
\[
f(AB) + f(AC) + f(BC) + f(AD) + f(BD) 
\ge  f(A) + f(B) + f(ABC) + f(ABD) + f(CD) 
\]
That is \emph{Ingleton's inequality}~\cite{Ing71}, 
a necessary condition for a
polymatroid to be linear.
It is a \emph{linear rank inequality},
that is, an inequality involving
a linear expression on the values of the rank function
that has to be satisfied by every linear polymatroid.

Given subsets $X, Y, Z$ of the ground set 
of a polymatroid $(E,f)$, we say that
$Z$ is a \emph{common information}
for the pair $(X,Y)$ if 
$f(Z) = f(X \lon Y)$ and
$f(Z | X) = f(Z | Y) = 0$.
A motivation for this concept and
its name is given in~\cite{HRSV00}.
A \emph{common information extension}
(or \emph{{\cip} extension} for short)
of $(E,f)$ for the pair $(X,Y)$ 
is an extension $(EZ,f)$
in which $Z$ is a commom information
for $(X,Y)$.
For every $\kos$-linear polymatroid
and every  pair of subsets of the ground set,
there exists 
a common information extension
that is also $\kos$-linear~\cite{HRSV00}.
Therefore, the existence of 
such extensions is another sufficient condition
for a polymatroid to be linear.

A stronger necessary condition is derived
by considering the existence of an arbitrary number
of iterated common information extensions.
Following the ideas in~\cite{Boll18} 
the concept of  \emph{$k$-{\cip} polymatroids}
was introduced in~\cite{BFP25}.
It is recursively defined as follows.
Every polymatroid is $0$-{\cip}.
For a positive integer $k$,  
a polymatroid is $k$-{\cip} if,
for every pair of subsets of the ground set, 
it admits a {\cip} extension that is a 
$(k-1)$-{\cip} polymatroid.
A polymatroid satisfies the
\emph{common information property},
or it is a \emph{{\cip} polymatroid},
if it is $k$-{\cip} for every 
positive integer $k$.

Every $1$-{\cip} polymatroid satisfies
the Ingleton's inequality~\cite{HRSV00}.
By using the common information property,
a number linear rank inequalities on variables 
were discovered~\cite{DFZ09}.
Since at most two {\cip} extensions
are needed to prove them, the inequalities
in~\cite{DFZ09}  are satisfied by every 
$2$-{\cip} polymatroid,
but some of them by every $1$-{\cip}
polymatroid too.

The common information property is 
closely related to other
necessary conditions for
a matroid to be linear, as the generalized 
Euclidean property
or the Levi's intersection property,
both discussed in~\cite{BaWa89,BFP25}.

Another kind of necessary conditions 
for a polymatroid to be linear
are derived from the tensor product
introduced in~\cite{Lov77,Mas77}.
A polymatroid $(E_1 \times E_2, g)$
is a \emph{tensor product}
of  the polymatroids $(E_1,f_1)$ and
$(E_2, f_2)$ if
$g(X \times Y) = f_1(X) f_2(Y)$
for every $X \subseteq E_1$ 
and $Y \subseteq E_2$. 
Not every pair of polymatroids admits a 
tensor product.
For example, there
is no tensor product of the
Vamos matroid and the uniform 
matroid $U_{2,3}$~\cite{LsV81}.
Nevertheless, by using the 
tensor product of vector spaces, 
a tensor product for every
pair of $\kos$-linear polymatroids 
can be found.
Therefore, the existence of a tensor 
product with $U_{2,3}$ is a necessary 
condition for a polymatroid to be linear.
Moreover,
it has been recently proved~\cite{BGILMS25}
that the Ingleton's inequality
is satisfied by every polymatroid that admits
a tensor product with $U_{2,3}$.
We present in the next section an
easier proof for that result. 
Nevertheless, it is superseded by
our main result (Theorem~\ref{st:main}).
Namely, every polymatroid that admits a tensor
product with $U_{2,3}$ is $1$-{\cip}.

\section{Yet another Proof of Ingleton's Inequality}

Consider a polymatroid $(E,f)$ and sets
$A_1, A_2, A_3 \subseteq E$.
Take $s = f(A_1 A_2 A_3)$ and  
$r_i = f(A_i)$ for $i = 1,2,3$.
In addition, for $\{i,j,k\} = \{1,2,3\}$, 
take
$s_i = f(A_j A_k)$ and  
$t_i =  r_j + r_k - s_i$.
Let $(V_x)_{x \in E}$ be a linear representation of
a polymatroid $(E,f)$.
Given three sets $A_1, A_2, A_3 \subseteq E$
consider, for each $i \in \{1,2,3\}$,
the vector subspace 
$U_i = \sum_{x \in A_i} V_x$.
An easy proof for the following well known inequalities
can be found, for example, in~\cite[Section~7.2]{FMP07}.
\begin{equation}
\label{eq:bas}
\max
\left\{
0, s - \sum s_i + \sum r_i 
\right\} \le
\dim(U_1 \cap U_2 \cap U_3) \le
\min\{t_1, t_2, t_3 \}
\end{equation}

\begin{proposition}
\label{st:reptens}
A linear representation $(V_x)_{x \in E}$ of $(E,f)$,
together with the obvious linear representation for
$U_{2,3}$, determine a tensor product 
$(E \times \{1,2,3\},g)$ of those
two polymatroids.
Consider
$A_1, A_2, A_3 \subseteq E$
and, for each $i \in \{1,2,3\}$,
the vector subspace 
$U_i = \sum_{x \in A_i} V_x$.
Then
\[
g(A_1^1 A_2^2 A_3^3) = 
f(A_1) + f(A_2) + f(A_3) - \dim(U_1 \cap U_2 \cap U_3)
\]
\end{proposition}

\begin{proof}
Use the properties of the tensor product of vector spaces.
\end{proof}

Take
\[
\alpha(A_1, A_2, A_3) = 
r_1 + r_2 + r_3 - \max 
\left\{
0, s - \sum s_i + \sum r_i 
\right\} = 
\min \left\{\sum r_i, \sum s_i - s \right\}
\]
and
\[
\beta(A_1, A_2, A_3) = 
r_1 + r_2 + r_3 - \min\{t_1, t_2, t_3 \} =
\max\{s_1 + r_1, s_2 + r_2, s_3 + r_3\}
\]
By combining the inequalities in~(\ref{eq:bas}) and
Proposition~\ref{st:reptens}, the bounds in~(\ref{bounds})
are satisfied under the hypotheses of Proposition~\ref{st:reptens}.
We prove next that those bounds apply to every tensor product
with $U_{2,3}$.

\begin{proposition}
\label{st:gentens}
Let $(E,f)$ be a polymatroid
that admits a tensor product
$(E \times \{1,2,3\},g)$ with $U_{2,3}$.
Then
\begin{equation}
\label{bounds}
\beta(A_1,A_2,A_3) \le  
g(A_1^1 A_2^2 A_3^3) \le
\alpha(A_1,A_2,A_3)
\end{equation}
for every $A_1,A_2,A_3 \subseteq E$.
\end{proposition}

\begin{proof}
Take $X = A_1^1 A_2^2 A_3^3$ and 
$Y = (A_1 A_2)^1 A_3^3$.
Then $X Y = (A_1 A_2)^1 A_2^2 A_3^3$ and
$X \cap Y = A_1^1 A_3^3$.
Observe that
$g(Y) = s_3 + r_3$ and
$g(X \cap Y) = r_1 + r_3$.
In addition
\[
g(X Y) = g ((A_1 A_2)^1 A_2^2 A_3^3) = 
g((A_1 A_2)^1 A_2^2 (A_2 A_3)^3) \ge
s_3 + s_1
\]
because $g(A_2^{1} A_2^{2}) = g(A_2^{1} A_2^{2} A_2^3)$.
Therefore, 
\[
g(X) + s_3 + r_3  
= g(X) + g(Y) 
\ge g(XY) + g(X \cap Y) 
\ge s_3 + s_1 + r_1+ r_3 
\]
and hence $g(X) \ge s_1 + r_1$.
By symmetry, the first inequality holds. 
Clearly,
\[
g(A_1^1 A_2^2 A_3^3) \le
g(A^1) + g(B^2) + g(C^3) = r_1 + r_2 + r_3
\]
Consider 
$Z = (A_1 A_3)^1 (A_2 A_3)^2  A_3^3$ and
$T = A_1^1 (A_2 A_3)^2 (A_1 A_2 A_3)^3$.
Then 
\[
ZT = (A_1 A_2 A_3 \times \{1,2,3\}) \sss A_2^1
\]
and $Z \cap T = X = A_1^1 A_2^2 A_3^3$.
Observe that $g(Z) = s_1 + s_2$
and $g(T) = s_3 + s$,
while $g(ZT) = 2 s$.
Therefore,
\[
s_1 + s_2 + s_3 + s
= g(Z) + g(T)
\ge g(ZT) + g(Z \cap T)
= 2 s + g(X) 
\]
which concludes the proof of the second inequality.
\end{proof}

\begin{proposition}
\label{st:newIng}
Every polymatroid
that admits a tensor product with $U_{2,3}$
satisfies the Ingleton's inequality.
\end{proposition}

\begin{proof}
Let $(E,f)$ be a polymatroid
that admits a tensor product
$(E \times \{1,2,3\},g)$ with $U_{2,3}$.
Let $A,B,C,D$ be four subsets of $E$.
Take $X = A^1 B^2 C^3$ and $Y = A^1 B^2 D^3$.
Then
\begin{eqnarray}
g(X) + g(Y) 
&{}\le{}& \alpha(A,B,C) + \alpha (A,B,D) \nonumber \\
&{}\le{}& 2 f(AB) + f(AC) + f(BC) + f(AD) + f(BD) - f(ABC) - f(ABD)
\nonumber
\end{eqnarray}
Since $XY =  A^1 B^2 (CD)^3$ and 
$X \cap Y = A^1 B^2$, 
\begin{eqnarray}
g(XY) + g(X \cap Y) 
&{}\ge{}& \beta(A,B,CD) + f(A) + f(B) \nonumber\\
&{}\ge{}& f(AB) + f(CD) + f(A) + f(B) \nonumber
\end{eqnarray}
Therefore,
\[
2 f(AB) + f(AC) + f(BC) + f(AD) + f(BD) - f(ABC) - f(ABD) 
\ge f(AB) + f(CD) + f(A) + f(B) 
\]
and hence
\[
f(AB) + f(AC) + f(BC) + f(AD) + f(BD) 
\ge  f(A) + f(B) + f(ABC) + f(ABD) + f(CD) 
\]
and the proof is concluded. 
\end{proof}

\section{Common Information from Tensor Products}

We begin by explaining the main 
idea behind the proof of Theorem~\ref{st:main}.
Assume that the polymatroid $(E,f)$ is
linearly represented by a 
collection $(V_x)_{x \in E}$ of 
subspaces of some vector space $V$.
For every $A \subseteq E$, 
consider the subspace
$U_A = \sum_{x \in A} V_x$.
On one hand, a common information extension
$(Ez,f)$ for a pair $(X,Y)$ is obtained
by adding the subspace
$V_z = U_X \cap U_Y$  
to the linear representation of $(E,f)$.
For every $A \subset E$, 
\begin{eqnarray}
f(A z) &{}={}& 
\dim (V_z + U_A) \nonumber \\
&{}={}& 
f(z) + f(A) - \dim(V_z \cap U_A)  \nonumber \\
&{}={}& 
f(z) + f(A) - \dim(U_X \cap U_Y \cap U_A) \label{eq:cieq}
\end{eqnarray}
On the other hand, 
the linear representation of 
$(E,f)$ determines a tensor
product with $U_{2,3}$.
By Proposition~\ref{st:reptens},
its rank function satisfies
\begin{equation}
\label{eq:tpeq}
g(X^1 Y^2 A^3) =
f(X) + f(Y) + f(A) - \dim(U_X \cap U_Y \cap U_A)
\end{equation}
Combining~(\ref{eq:cieq}) and~(\ref{eq:tpeq}),
\[
f(Az) = f(z) + f(A) + g(X^1 Y^2 A^3) - f(X) - f(Y) - f(A)
\] 
and, by taking into account that 
$f(z) = f(X) + f(Y) - f(XY)$,
\[
f(Az) =  g(X^1 Y^2 A^3) - f(XY)
\] 
That equality tells how to find common information
extensions from a tensor product with $U_{2,3}$.

\begin{theorem}
\label{st:main}
Every polymatroid that admits a tensor product
with $U_{2,3}$ is $1$-\cip.
\end{theorem}

\begin{proof}
Let $(E,f)$ be a polymatroid that admits
a tensor product $(E \times \{1,2,3\}, g)$
with $U_{2,3}$.
Given a pair $(X,Y)$ of subsets of $E$
take the extension $(Ez,f)$ of $(E,f)$ determined by
\[
f(Az) = g(X^1 Y^2 A^3) - f(XY)
\]
for every $A \subseteq E$.
We prove next that $(Ez,g)$ is a polymatroid
in which $z$ is a common information for the pair $(X,Y)$.
We prove first that
$f(z) = f( X\lon Y)$ and $f(z| X) = f(z| Y) = 0$.
Indeed, 
\[
f(z) = g(X^1 Y^2) - f(XY)
= f(X) + f(Y) - f(XY) = f(X \lon Y)
\]
and
\[
f(Xz) = g(X^1 Y^2 X^3) - f(XY) = f(X) + f(XY) - f(XY) = f(X)
\]
which implies that $f(z|X) = 0$.
Symmetrically, $f(z| Y) = 0$.
In order to prove that $(Ez, f)$ is a polymatroid, we 
check that $f(x \lon y |A) \ge 0$
for all $A \subseteq E z$ and $x,y \in E z \sss A$.
Since $(E,f)$ is a polymatroid, that holds if
$z \notin A x y$. 
If $A = Bz$ with $B \subseteq E$, then
\begin{eqnarray}
f(x \lon y | A) 
&{}={}&
f(x \lon y | Bz) \nonumber \\
&{}={}&
f(Bxz) + f(Byz) - f(Bxyz) - f(Bz)  \nonumber \\
&{}={}&
g(X^1 Y^2 (Bx)^3) +
g(X^1 Y^2 (By)^3) -
g(X^1 Y^2 (Bxy)^3) -
g(X^1 Y^2 B^3) 
\nonumber \\
&{}\ge{}& 0 \nonumber
\end{eqnarray}
where the last inequality follows
from the submodularity of $g$. 
If $x \ne z$,
\begin{eqnarray}
f(x \lon z | A) 
&{}={}&
f(Ax) + f(Az) - f(Axz) - f(A)  \nonumber \\
&{}={}&
f(Ax) + g(X^1 Y^2 A^3) - 
g(X^1 Y^2 (Ax)^3) - f(A) \nonumber \\
&{}={}&
g((Ax)^3) +
g(X^1 Y^2 A^3) -
g(X^1 Y^2 (Ax)^3) -
g(A^3)\nonumber \\
&{}\ge{}& 0 \nonumber
\end{eqnarray}
because $g$ is submodular.
Finally,
\[
f(z \lon z |A) = f(Az) - f(A) = 
g(X^1 Y^2 A^3) - f(XY) - f(A) \ge
g(X^1 Y^2 A^3) - \alpha(X,Y,A) \ge 0
\]
and the proof is concluded.
\end{proof}

\end{document}